\documentclass[11pt]{article}
\usepackage[margin=1in]{geometry}
\usepackage{lmodern}
\usepackage{microtype}

\usepackage{amsmath,amssymb,amsthm}
\usepackage{mathtools}

\usepackage{algorithm}
\usepackage{algpseudocode}

\usepackage{graphicx}
\graphicspath{{figures/}}
\usepackage{booktabs}
\usepackage{array}
\usepackage[font=small,labelfont=bf]{caption}
\usepackage{xcolor}
\usepackage{subcaption,comment}

\usepackage{pgfplots}
\usepgfplotslibrary{groupplots}
\pgfplotsset{compat=1.18}

\usepackage[noblocks]{authblk}

\usepackage[numbers,sort&compress]{natbib}
\usepackage[hidelinks]{hyperref}
\usepackage[nameinlink,capitalise,noabbrev]{cleveref}

\crefformat{equation}{#2(#1)#3}
\Crefformat{equation}{#2(#1)#3}
\crefrangeformat{equation}{#3(#1)#4--#5(#2)#6}
\Crefrangeformat{equation}{#3(#1)#4--#5(#2)#6}
\crefmultiformat{equation}{#2(#1)#3}{ and~#2(#1)#3}{, #2(#1)#3}{ and~#2(#1)#3}
\Crefmultiformat{equation}{#2(#1)#3}{ and~#2(#1)#3}{, #2(#1)#3}{ and~#2(#1)#3}

\makeatletter
\providecommand{\theHALG@line}{\thealgorithm.\arabic{ALG@line}}
\renewcommand{\theHALG@line}{\thealgorithm.\arabic{ALG@line}}
\makeatother

\theoremstyle{plain}
\newtheorem{theorem}{Theorem}
\newtheorem{lemma}{Lemma}
\newtheorem{proposition}{Proposition}
\newtheorem{corollary}{Corollary}

\theoremstyle{definition}
\newtheorem{definition}{Definition}

\theoremstyle{remark}
\newtheorem{remark}{Remark}

\newcommand{\DCSG}{\texttt{DC-SG}}
\newcommand{\R}{\mathbb{R}}
\DeclareMathOperator{\dist}{dist}

\DeclareMathOperator*{\argmin}{arg\,min}

\usepackage{titlesec}

\begin{document}

\title{
\bfseries\Large
Distance-Certified Polyak-type Stepsizes for Subgradient Methods under H{\"o}lder Growth\thanks{This work is supported by NSF CAREER 2340713. The work of K. Weraarchakul is also supported by the ECE Spring/Summer Research Internship Program (SRIP) at UC San Diego. Emails: \texttt{\{kweraarchakul,zhengy\}@ucsd.edu}.}
}
\author[1]{Kornpholkrit Weraarchakul}
\author[1]{Yang Zheng}
\affil[1]{\small Department of Electrical and Computer Engineering, University of California San Diego}
\date{\small \today \vspace{-5ex}}

\maketitle

\begin{abstract}
    In this paper, we consider nonsmooth convex optimization under a H{\"o}lder growth condition, with quadratic growth and weak sharp minima as two important special cases. Existing subgradient methods often employ predetermined decaying stepsize schedules, whereas the classical Polyak stepsize selects its stepsize adaptively but requires knowledge of the optimal value. We introduce a new distance-certified subgradient method (\DCSG) that follows a one-step optimization principle analogous to that underlying the Polyak stepsize, without requiring the optimal value. Specifically, \DCSG{} maintains a valid upper certificate for the current squared distance to the solution set. The method selects each stepsize by minimizing the worst-case upper bound on the next squared distance over all distances consistent with this certificate. We show that the resulting minimax problem admits a closed-form solution, yielding an explicit adaptive stepsize update. Under uniformly bounded subgradients, \DCSG{} improves upon the convergence rate of the classical subgradient method and matches the state-of-the-art rates for Lipschitz convex optimization under H{\"o}lder growth. In particular, it achieves linear convergence under sharp growth. Numerical experiments illustrate the performance of \DCSG.
\end{abstract}

\section{Introduction}
We consider the convex optimization problem
\begin{equation}\label{eq:convex-optimization}
    \min_{x\in C}\; f(x),
\end{equation}
where $f:\R^n\to\R \cup\, \{+\infty\}$ is a closed convex function and $C\subseteq\mathrm{dom}\,f$ is a closed and convex set. We do not assume that $f$ is smooth or strongly convex.
Let $X^\star:=\argmin_{x\in C}f(x)$ denote the solution set and $f^\star:=\min_{x\in C}f(x)$ the optimal value. Throughout this paper, we assume that $X^\star$ is nonempty and $f^\star$ is finite. Problem \cref{eq:convex-optimization} arises in many applications such as machine learning, compressed sensing, image processing, computer vision, and many others \cite{Nesterov2004,BoydVandenberghe2004,luo2006introduction}.

When $f$ is smooth, projected gradient descent and its accelerated variants \cite{Nesterov2004} have been well-developed to solve \cref{eq:convex-optimization}. For nonsmooth objectives, the subgradient method is among the simplest first-order algorithms. Its development dates back to foundational work in the 1960s \cite{Shor1968,Polyak1969}, and it remains widely used because of its simple implementation and low per-iteration cost. When $f$ is globally Lipschitz-continuous, the classical subgradient method achieves a convergence rate of $f(\hat{x}_k) - f^\star \leq O(1/\sqrt{k})$ after $k$ subgradient evaluations for a suitably averaged point $\hat{x}_k$ \cite{Nesterov2004}. 
This rate is known to be the worst-case optimal for general nonsmooth convex optimization among the first-order oracle model \cite{NemirovskiYudin1983,Nesterov2004}. Faster rates are nevertheless possible for certain classes of $f$ with additional structural information   \cite{JohnstoneMoulin2020,goffin1977convergence,YangLin2018,Wang2025,Polyak1969,davis2025gradient,liao2026error,liao2025proximal,liao2026revisiting}.

In this paper, we consider Lipschitz convex optimization under the H{\"o}lder growth condition
\begin{equation} \label{eq:HEB}
    f(x) - f^\star \geq c\,\mathrm{dist}(x,X^\star)^p, \quad \forall x \in C,
\end{equation}
for some choice of $c>0$ and $p\geq 1$. This is equivalent to
$
    \dist(x,X^\star)
    \leq
    c^{-1/p}\bigl(f(x)-f^\star\bigr)^{1/p}, \forall x \in C.
$ 
This equivalent form bounds the distance to the solution set in terms of the objective gap and is commonly referred to as a H{\"o}lder error bound. The error-bound exponent is often denoted by $\theta\in(0,1]$; under this parameterization, we have $p=1/\theta$ as used in \cite{JohnstoneMoulin2020,Wang2025}. This  H{\"o}lder growth property is also called the \L{}ojasiewicz inequality \cite{bolte2007lojasiewicz} and is closely related to the Kurdyka-\L{}ojasiewicz (KL) inequality \cite{BolteNguyenPeypouquetSuter2017}.  H{\"o}lder error bounds have been widely studied in optimization, sensitivity analysis, and convergence-rate analysis
\cite{JohnstoneMoulin2020,BolteNguyenPeypouquetSuter2017,zhou2017unified,bolte2007lojasiewicz,liao2026error,
BurkeFerris1993,luo1993error}. 

The H{\"o}lder growth condition holds for many problems arising in applications. In particular, every proper closed semialgebraic convex function satisfies a local error bound of the form \cref{eq:HEB} for suitable constants $c>0$ and $p\geq1$ \cite{BolteNguyenPeypouquetSuter2017}. This additional growth structure enables many first-order algorithms to achieve faster convergence rates \cite{JohnstoneMoulin2020,goffin1977convergence,YangLin2018,Wang2025,
Polyak1969,davis2025gradient,liao2026error,liao2025proximal}. Two particularly important cases are $p=1$ and $p=2$. The case $p=2$ is commonly known as the quadratic growth (QG) condition \cite{zhou2017unified,liao2026error}, and the case $p=1$ corresponds to the sharp growth condition, also known as weak sharp minima \cite{BurkeFerris1993,JohnstoneMoulin2020}. The QG condition has been used to establish linear convergence of proximal point and related methods \cite{liao2026error}. Moreover, for convex functions, QG is equivalent to several other commonly used regularity conditions, including the restricted secant inequality, the subdifferential error bound, and the Polyak--{\L}ojasiewicz condition
\cite[Theorem~3.1]{liao2026error}.

Under the H{\"o}lder growth condition \cref{eq:HEB}, several techniques have been developed to accelerate the convergence of subgradient methods across a range of settings. Representative approaches~include predetermined decaying stepsize schedules \cite{JohnstoneMoulin2020,goffin1977convergence,Wang2025}, restarting schemes \cite{RouletAspremont2020,YangLin2018} and Polyak-type stepsizes \cite{Polyak1969,DevanathanBoyd2024,he2025new,zamani2025exact}. These approaches substantially improve convergence rates of subgradient methods by carefully designing the stepsizes. Nevertheless, the predetermined schedules used by restarting and decaying schemes \cite{JohnstoneMoulin2020,goffin1977convergence,Wang2025,RouletAspremont2020,YangLin2018} present several practical and conceptual limitations. Their stepsizes are typically tuned using global problem information, such as a uniform subgradient bound and growth parameters. To achieve the fast rates, it may require carefully chosen stage lengths and different decaying coefficients. This makes the resulting schedules relatively intricate and obscures the connection between each stepsize and the corresponding one-step progress of the subgradient method. Furthermore, these schemes commonly use a global Lipschitz bound for stepsize tuning. This choice may lead to conservative stepsize choices, and consequently, these predetermined schemes \cite{JohnstoneMoulin2020,Wang2025,YangLin2018} do~not capture the additional potential acceleration when the subgradients diminish near the solution~set. 

The classical Polyak stepsize \cite{Polyak1969} follows a different principle and does not rely on a predetermined schedule. At each iteration, it minimizes a one-step upper bound on the squared distance
to the solution set using the current objective gap and subgradient norm. This allows the stepsize to adapt naturally to the local problem geometry and may permit more aggressive steps than standard constant or
diminishing stepsizes \cite{zamani2025exact,he2025new}. Recently, Wang \cite{Wang2025} showed that the classical Polyak stepsize attains the minimax-optimal convergence rate for reducing the distance to the
solution set over globally Lipschitz convex functions satisfying \cref{eq:HEB}. Implementing the Polyak stepsize, however, requires exact knowledge of the optimal value $f^\star$, which may be unavailable. 
Note that avoiding exact knowledge of the optimal value requires special treatment; see, for example, dynamic target-level schemes \cite{GoffinKiwiel1999,NedicBertsekas2001,LiuBraginChenGuan2025} and restarted procedures that successively refine a lower bound on $f^\star$ \cite{HazanKakade2019}.

In this paper, we introduce a new distance-certified subgradient method (\DCSG) for Lipschitz convex optimization under  H{\"o}lder growth. Instead of prescribing a decaying stepsize schedule in advance, \DCSG{} follows a one-step optimization principle analogous to that underlying the classical Polyak stepsize. Specifically, our key idea is to maintain a valid upper certificate for the current squared distance to the solution set. We select the stepsize by minimizing the worst-case upper bound on the next squared distance over all distances consistent with the current certificate. This construction uses the H{\"o}lder growth parameters, but avoids the requirement of the optimal value $f^\star$. We show that the resulting certificates remain valid and nonincreasing (\Cref{prop:certificate-monotonicity}) and that the minimax principle admits an explicit adaptive solution determined by the current certificate, the observed subgradient norm, and the  H{\"o}lder growth parameters (\Cref{thm:dc-sg-certificate-1-leq-p-leq-2,thm:dc-sg-certificate-p-geq-2}). 

Similar to the Polyak stepsize, \DCSG{} adapts its stepsize at every iteration and requires no predetermined schedule.
Also, the implementation of  \DCSG{} requires no global Lipschitz constant, which enters only the worst-case convergence analysis. Under uniformly bounded subgradients, we prove that \DCSG{} provides the \textit{last-iterate} guarantee
$$
    \dist(x_k,X^\star)
    \leq O\left(k^{-\frac{1}{2(p-1)}}\right)
$$
when \(p>1\) (\cref{thm:polynomial-certificate-decay-lipschitz}) and converges linearly when \(p=1\) (\cref{thm:geometric-certificate-decay-p-equals-1}). These rates are faster than the classical subgradient guarantee and remove the additional logarithmic factor appearing in the restarted subgradient method \cite{YangLin2018}. Also, they match the state-of-the-art rates for Lipschitz convex optimization under H{\"o}lder growth \cite{JohnstoneMoulin2020,Wang2025}. Finally,  \DCSG{} automatically achieves faster rates when the norm of observed subgradients diminishes near the solution set (\cref{cor:shrinking-subgradient-rates}). This~additional acceleration is not captured by the predetermined decaying stepsize schedules \cite{JohnstoneMoulin2020,Wang2025}. We present further comparisons with the state-of-the-art in \Cref{subsection:comparison}.

The remainder of this paper is organized as follows. \Cref{sec:notations-and-preliminaries} introduces the necessary notation and preliminary results. \Cref{sec:distance-certified-subgradient-method} develops the distance-certified subgradient method (\DCSG). \Cref{sec:dc-sg-shrinking-rates} establishes its convergence rates and compares them with state-of-the-art guarantees. \Cref{section:experiments} illustrates the numerical performance of \DCSG. We conclude the paper in \Cref{section:conclusion}. 

\section{Preliminaries}\label{sec:notations-and-preliminaries}

Recall that the subdifferential of $f$ at $x$ is defined by
$$
    \partial f(x)=\{g\in\mathbb{R}^n \mid f(y)\geq f(x)+\langle g,y-x\rangle,\ \forall y\in\mathbb{R}^n\}.
$$
Given an initial point $x_0\in C$, the projected subgradient method generates iterates according to 
$$
    x_{k+1}=P_C(x_k-\alpha_k g_k),\qquad k\geq0,
$$
where $g_k\in\partial f(x_k)$, $\alpha_k>0$ is the stepsize, and $P_C$ denotes the Euclidean projection onto $C$. 
The method is listed in \Cref{alg:projected-subgradient-method}. Different stepsize choices lead
to different convergence guarantees \cite{JohnstoneMoulin2020,goffin1977convergence,YangLin2018,Wang2025,Polyak1969,davis2025gradient}. Among them, the Polyak stepsize \cite{Polyak1969} is particularly relevant to our development.

\begin{algorithm}[t]
    \caption{Projected subgradient method}
    \label{alg:projected-subgradient-method}
    \begin{algorithmic}[1]
        \Require An initial point $x_0\in C$, stepsizes $\alpha_k>0$, and an iteration limit $N$.
        \For{$k=0,1,\ldots,N-1$}
            \State Choose $g_k\in\partial f(x_k)$.
            \State Set $x_{k+1}=P_C(x_k-\alpha_k g_k)$.
        \EndFor
    \end{algorithmic}
\end{algorithm}

We focus on functions satisfying the following H{\"o}lder growth condition.
\begin{definition}[H{\"o}lder growth]\label{def:holderian-error-bound}
    A function $f$ satisfies $(p,c)$-H{\"o}lder growth on $C$ if there exist $c>0$ and $p\geq1$ such that
    \begin{equation}\label{eq:holderian-error-bound}
        f(x)-f^\star\geq c\,\dist(x,X^\star)^p,\qquad \forall x\in C,
    \end{equation}
    where $\dist(x,X^\star):=\inf_{y\in X^\star}\|x-y\|$ denotes the distance from $x$ to the solution set.
\end{definition}

We refer to \cref{eq:holderian-error-bound} as the H{\"o}lder growth condition.
This condition is also known as a H{\"o}lderian error bound and is closely related to the Kurdyka-\L{}ojasiewicz (KL) inequality \cite{BolteNguyenPeypouquetSuter2017}.
 
It is well-known that the projected subgradient method may not decrease the objective value at every iteration. Still, its progress can be measured by the distance of the iterates to the solution set.
For each iteration of \Cref{alg:projected-subgradient-method}, we define
\begin{equation} \label{eq:distance-def}
    e_k:=\dist(x_k,X^\star)^2,     \qquad     h_k:=\|g_k\|,
\end{equation}
which denote the squared-distance error and the norm of the selected (or observed) subgradient, respectively. The following standard lemma provides the fundamental one-step improvement, which underlies most convergence analyses of subgradient methods \cite{Nesterov2004,JohnstoneMoulin2020,goffin1977convergence,YangLin2018,Wang2025,Polyak1969,davis2025gradient}. We include a brief proof to highlight its simplicity and also for completeness.

\begin{subequations}
\begin{lemma}[Squared-distance recursion]\label{lem:distance-recursion}
    Let $f$ be a closed convex function. For every iterate $x_k$ generated by \Cref{alg:projected-subgradient-method}, we have
    \begin{equation} \label{eq:distance-function-value}
        e_{k+1}
        \leq e_k
        -2\alpha_k\bigl(f(x_k)-f^\star\bigr)
        +\alpha_k^2h_k^2.
    \end{equation}
    If, in addition, $f$ satisfies the $(p,c)$-H{\"o}lder growth condition~\cref{eq:holderian-error-bound}, then
    \begin{equation}\label{eq:distance-recursion}
        e_{k+1}\leq e_k-2\alpha_kce_k^{p/2}+\alpha_k^2h_k^2. 
    \end{equation}
\end{lemma}
\end{subequations}

\begin{proof}
    Fix any $x^\star\in X^\star$. Since $x^\star\in C$, we have $P_C(x^\star)=x^\star$. By the non-expansiveness of the projection operator, we have 
    $$
    \begin{aligned}
        \|x_{k+1}-x^\star\|^2
        &=\|P_C(x_k-\alpha_k g_k)-P_C(x^\star)\|^2\\
        &\leq\|x_k-\alpha_k g_k-x^\star\|^2\\
        &=\|x_k-x^\star\|^2-2\alpha_k\langle g_k,x_k-x^\star\rangle+\alpha_k^2\|g_k\|^2.
    \end{aligned}
    $$
    Since $g_k\in\partial f(x_k)$, the subgradient inequality gives $f(x_k)-f^\star\leq\langle g_k,x_k-x^\star\rangle$. Recalling the definition that $h_k=\|g_k\|$ in \cref{eq:distance-def}, we obtain
    $$
        \|x_{k+1}-x^\star\|^2
        \leq
        \|x_k-x^\star\|^2-2\alpha_k(f(x_k)-f^\star)+\alpha_k^2h_k^2.
    $$
    Taking the infimum over $x^\star\in X^\star$ on both sides gives \cref{eq:distance-function-value}.
    If $f$ satisfies the $(p,c)$-H{\"o}lderian~condition~\cref{eq:holderian-error-bound}, then
    $
        f(x_k)-f^\star\geq c\,\dist(x_k,X^\star)^p=ce_k^{p/2}.
    $
    Substituting this bound into \cref{eq:distance-function-value} gives~\cref{eq:distance-recursion}.
\end{proof}

From \cref{eq:distance-function-value,eq:distance-recursion}, the negative term scales linearly with $\alpha_k$, while the positive term scales quadratically.
Thus, whenever $x_k$ is nonoptimal, a sufficiently small stepsize guarantees $e_{k+1}<e_k$.
Building on this insight, constant stepsizes and predetermined decaying stepsizes can yield complexity guarantees for finding an iterate $x_k$ satisfying $e_k \leq \epsilon$ for functions with H{\"o}lder growth \cite{JohnstoneMoulin2020,Wang2025}. 

\begin{remark} \label{remark:upper-bound}
If $f$ is globally $L$-Lipschitz continuous, then $h_k\leq L$ for every iteration. Consequently, \cref{eq:distance-recursion} reduces to
    $
        e_{k+1} \leq e_k-2\alpha_k c e_k^{p/2}+\alpha_k^2L^2,
    $ 
    which has been widely used in existing stepsize designs \cite{JohnstoneMoulin2020,Wang2025}. However, this worst-case bound does not exploit the subgradient norm $h_k$ observed at each iteration. As a result, the algorithms in \cite{JohnstoneMoulin2020,Wang2025} may fail to adapt to smaller subgradients and may not achieve faster convergences. We provide further discussion in \Cref{remark:norm-aware-stepsizes} later. \hfill $\square$
\end{remark}

\section{A distance-certified subgradient method (\DCSG)}\label{sec:distance-certified-subgradient-method}

In this section, we develop a new distance-certified subgradient method for convex functions satisfying H{\"o}lder growth. We first introduce the underlying stepsize design, then derive explicit stepsize updates, and finally provide technical proofs.

\subsection{Distance-certified stepsize design} \label{subsection:stepsize-design}
The stepsize design plays a fundamental role in the convergence of the projected subgradient method in \Cref{alg:projected-subgradient-method}. Indeed, once the oracle selects a subgradient $g_k$, the stepsize $\alpha_k$ is the only remaining design choice and fully determines the resulting update. Common choices include constant stepsizes $\alpha_k=\alpha$ \cite{Nesterov2004}, diminishing stepsizes satisfying
$
    \sum_{k=0}^{\infty}\alpha_k=\infty, \sum_{k=0}^{\infty}\alpha_k^2<\infty,
$ 
such as $\alpha_k=a/(b+k)$ with $a,b>0$ \cite{Nesterov2004,BoydVandenberghe2004}, and geometrically decaying stepsizes \cite{JohnstoneMoulin2020,goffin1977convergence,Wang2025}. 

\begin{subequations}
If the optimal value $f^\star$ is known, minimizing the right-hand side of the squared-distance recursion \cref{eq:distance-function-value} with respect to $\alpha_k$ yields the celebrated Polyak stepsize \cite{Polyak1969}
\begin{equation} \label{eq:polyak-stepsize}
\alpha_{\texttt{Polyak},k} = \frac{f(x_{k}) - f^\star}{h_k^2},\quad k = 0, 1, 2, \ldots.
\end{equation}
This stepsize can be interpreted as a greedy choice that minimizes the one-step upper bound on the squared-distance error. By incorporating both the current optimality gap and the observed subgradient norm, it naturally adapts to the local optimization landscape and can be more aggressive than standard constant or diminishing stepsizes \cite{zamani2025exact,he2025new}. 
Following the same idea as the Polyak stepsize, suppose that $f$ satisfies the $(p,c)$-H{\"o}lder growth condition and that the exact squared-distance error $e_k$ is available. Minimizing the right-hand side of \cref{eq:distance-recursion} with respect to $\alpha_k$~yields
\begin{equation} \label{eq:distance-based_Polyak-type_stepsize}
    \alpha_{\texttt{dist},k}=\frac{ce_k^{p/2}}{h_k^2},\quad k = 0, 1, 2, \ldots. 
\end{equation}
We refer to this choice as an exact distance-based Polyak-type stepsize. A related idea was recently explored in \cite{Wang2025}.
\end{subequations}

In general, the exact squared-distance error $e_k$ is unavailable, as its computation typically~requires knowledge of the solution set $X^\star$. Our key idea is to replace $e_k$ with a certified upper bound $E_k$
satisfying $e_k\leq E_k$ and choose the stepsize that minimizes the worst-case one-step upper bound in \cref{eq:distance-recursion} over all $e_k\in[0,E_k]$. Specifically, for $E\geq 0$, define
\begin{equation}\label{eq:one-step-worst-case-certificate}
    \Phi_k(\alpha,E):=\alpha^2h_k^2+\max_{0\leq e\leq E}\left\{e-2\alpha ce^{p/2}\right\},
\end{equation}
which represents the worst-case one-step upper bound in \cref{eq:distance-recursion} over all $e_k\in[0,E]$.  We then select
\begin{equation}\label{eq:stepsize-optimization}
    \alpha_{\texttt{cert},k} = \argmin_{\alpha\geq0}\Phi_k(\alpha,E_k), \quad k = 0, 1, 2, \ldots, 
\end{equation}
which we refer to as a distance-certified Polyak-type stepsize. By
construction, the value 
\begin{equation} \label{eq:stepsize-upper-bound}
    U_{k+1} := \Phi_k(\alpha_{\texttt{cert},k},E_k)
\end{equation}
provides an upper bound on the squared-distance error of the next iterate
$
    x_{k+1} = P_C\bigl(x_k-\alpha_{\texttt{cert},k}g_k\bigr),
$ 
namely, $e_{k+1}\leq U_{k+1}$. We may thus choose the next certified upper bound $E_{k+1} = U_{k+1}$ and repeat \Cref{eq:stepsize-optimization,eq:stepsize-upper-bound}. Starting from an initial bound $U_0\geq e_0$, this generates a sequence of certified bounds satisfying $e_k\leq E_k$ for all $k\geq 0$ and forms the basis of our distance-certified subgradient~method.

In addition, the observed subgradient $g_k$ can provide a further upper bound on the squared-distance error $e_k$. This observation is formalized in the following lemma.

\begin{lemma}[Subgradient-based distance certificate]
\label{lem:subgradient-distance-certificate}
    Suppose that $f$ is convex and satisfies the $(p,c)$-H{\"o}lder growth condition~\cref{eq:holderian-error-bound}. For any $x\in C$ and $g\in\partial f(x)$, we have
    $
        c\,\dist(x,X^\star)^p \leq \|g\|\,\dist(x,X^\star).
    $ 
    Consequently, if $p>1$, then
    \begin{equation}\label{eq:subgradient-distance-certificate}
        \dist(x,X^\star)^2 \leq \left(\frac{\|g\|}{c}\right)^{\frac{2}{p-1}}.
    \end{equation}
    If $p=1$ and $x\notin X^\star$, then $\|g\|\geq c$.
\end{lemma}

\begin{proof}
    Let $x^\star=P_{X^\star}(x)$, so that $\|x-x^\star\|=\dist(x,X^\star)$. By the subgradient inequality
    and the Cauchy--Schwarz inequality, we have $
        f(x)-f^\star \leq \langle g,x-x^\star\rangle \leq
        \|g\|\,\dist(x,X^\star).
    $ 
    Combining this inequality with \cref{eq:holderian-error-bound} proves the first statement. 

    If $x\in X^\star$, \Cref{eq:subgradient-distance-certificate} holds trivially. Otherwise, we have
    $
        c\,\dist(x,X^\star)^{p-1}\leq\|g\|.
    $ 
    For $p>1$, rearranging and squaring this inequality establishes \cref{eq:subgradient-distance-certificate}. For $p=1$, it reduces to $\|g\|\geq c$.
\end{proof}

\begin{algorithm}[t]
    \caption{Distance-certified subgradient method (\DCSG)}
    \label{alg:dc-sg}    
    \begin{algorithmic}[1]
        \Require An initial point $x_0\in C$, an initial bound $U_0\geq e_0$, parameters $c>0,\,p\geq 1$, and an iteration limit $N$.
        \For {$k= 0,1,\ldots,N-1$}
            \State Choose $g_k\in\partial f(x_k)$ and set $h_k=\|g_k\|$; 
            \State Compute $E_k$ according to  \cref{eq:certificate-update};
            \State Compute $\alpha_{\texttt{cert},k}$ according to \cref{eq:stepsize-optimization};
            \State Set $U_{k+1}$ according to \Cref{eq:stepsize-upper-bound};
            \State Set $x_{k+1}=P_C(x_k-\alpha_{\texttt{cert},k}g_k)$;
        \EndFor
    \end{algorithmic}
\end{algorithm}

After observing the selected subgradient $g_k$, \Cref{lem:subgradient-distance-certificate} allows us to
refine the upper bound $U_k$ from \cref{eq:stepsize-upper-bound} by setting
\begin{equation}\label{eq:certificate-update}
    E_k := \begin{cases}
        \min\left\{
        U_k,\,
        \left({h_k}/{c}\right)^{\frac{2}{p-1}}
    \right\} & \text{if} \quad p > 1,\\
    U_k, & \text{if} \quad p = 1. 
    \end{cases}
\end{equation}
By construction, we always have $e_k\leq E_k\leq U_k$. Combining the certificate update \cref{eq:certificate-update} with the stepsize design \cref{eq:stepsize-optimization} and the bound propagation \cref{eq:stepsize-upper-bound}, we present the distance-certified subgradient method (\DCSG), listed in \Cref{alg:dc-sg}.
We have a simple observation for this algorithm. 
\begin{proposition}[Monotonicity of distance certificates]
\label{prop:certificate-monotonicity}
    Suppose that $U_0\geq e_0$. Then the iterates generated by \Cref{alg:dc-sg} satisfy
    \[
        e_{k+1} \leq E_{k+1} \leq U_{k+1} \leq E_k \leq U_k,  \qquad \forall k\geq 0.
    \]
    In particular, both $\{E_k\}$ and $\{U_k\}$ are nonincreasing.
\end{proposition}

\begin{proof}
    By the definition of $\alpha_{\texttt{cert},k}$ in \cref{eq:stepsize-optimization}, we have
    \[
        U_{k+1} = \min_{\alpha\geq 0}\Phi_k(\alpha,E_k) \leq \Phi_k(0,E_k) = E_k.
    \]
    The remaining inequalities follow directly from \Cref{lem:distance-recursion,lem:subgradient-distance-certificate} and the certificate update \cref{eq:certificate-update}.
\end{proof}

In \Cref{sec:dc-sg-shrinking-rates}, we show that \DCSG{} achieves state-of-the-art convergence rates for solving \cref{eq:convex-optimization} with convex functions under the H{\"o}lder growth condition \cref{eq:holderian-error-bound}, while offering additional benefits.

\subsection{Explicit distance-certified stepsize updates}

The distance-certified stepsize $\alpha_{\texttt{cert},k}$ is defined by minimizing the worst-case upper bound in \cref{eq:stepsize-optimization}. Fortunately, the scalar structure of this minimax problem allows an explicit characterization of its solution. We first observe that the stepsize optimization problem \cref{eq:one-step-worst-case-certificate} is strongly convex and thus admits a unique minimizer.

\begin{lemma}[Unique distance-certified stepsize] 
    Given any $E_k>0$, define the inner worst-case function 
    \begin{equation} \label{eq:maximum-function}
        M_k(\alpha):=\max_{0\leq e\leq E_k}\left\{e-2\alpha c e^{p/2}\right\} 
    \end{equation}
    where $\alpha\geq0, p \geq 1.$ Then $M_k$ is convex and nonnegative on $[0,\infty)$. Consequently, for any $h_k>0$, 
    $
        \alpha\mapsto\Phi_k(\alpha,E_k)
        =
        \alpha^2h_k^2+M_k(\alpha)
    $
    is strongly convex in $\alpha$ and admits a unique minimizer on $[0,\infty)$.
\end{lemma}

\begin{proof}
    For each fixed $e\in[0,E_k]$, the function $\alpha\mapsto e-2\alpha c e^{p/2}$ is affine. Therefore, $M_k$ is convex as the pointwise maximum of affine functions. Moreover, since $e=0$ is feasible in the inner maximization, we have $M_k(\alpha)\geq 0$ for all $\alpha\geq 0$. 
    
    If $h_k>0$, then $\alpha\mapsto\alpha^2h_k^2$ is strongly convex. Hence, its sum with the convex function $M_k$ is also strongly convex. Furthermore,
    $
        \Phi_k(\alpha,E_k)  \geq \alpha^2h_k^2  \to\infty
        \, \text{as }\alpha\to\infty.
    $ 
    Therefore, $\Phi_k(\cdot,E_k)$ attains its minimum on $[0,\infty)$, and strong convexity guarantees that the minimizer is unique.
\end{proof}

\begin{figure}[t]
  \centering
  \resizebox{0.9\linewidth}{!}{
\begin{tikzpicture}
  \definecolor{phiEndpoint}{HTML}{E53935}
  \definecolor{phiSecond}{HTML}{1E88E5}
  \begin{groupplot}[
    group style={
      group size=2 by 1,
      horizontal sep=0.9cm,
      ylabels at=edge left,
      yticklabels at=edge left
    },
    width=7.6cm,
    height=4.2cm,
    scale only axis,
    xmin=0, xmax=1.25,
    ymin=0, ymax=1.8,
    xtick={0,0.4,0.8,1.2},
    ytick={0,0.5,1,1.5},
    axis lines=left,
    axis line style={line width=0.6pt,-},
    tick align=outside,
    tick style={line width=0.6pt},
    tick label style={font=\small},
    title style={font=\normalsize,yshift=2pt},
    label style={font=\normalsize},
    xlabel={Stepsize $\alpha$},
    ylabel={$\Phi_k(\alpha,E_k)$},
    ymajorgrids=true,
    grid style={gray!18,line width=0.3pt},
    samples=201,
    clip=true,
    /pgf/number format/fixed,
    /pgf/number format/precision=1,
    /pgf/number format/fixed zerofill,
    endpoint valid/.style={phiEndpoint,line width=1.2pt},
    endpoint extension/.style={phiEndpoint,densely dashed,line width=0.85pt},
    second valid/.style={phiSecond,line width=1.2pt},
    second extension/.style={phiSecond,densely dashed,line width=0.85pt}
  ]

    \nextgroupplot[]
    \addplot[endpoint extension,domain=0.5:1.25] {(x-1)^2};
    \addplot[second extension,domain=0:0.5] {x^2};
    \addplot[endpoint valid,domain=0:0.5] {(x-1)^2};
    \addplot[second valid,domain=0.5:1.25] {x^2};

    \nextgroupplot[]
    \addplot[endpoint extension,domain=0.25:1.25] {(x-1)^2};
    \addplot[second extension,domain=0.02:0.25] {x^2+1/(8*x)};
    \addplot[endpoint valid,domain=0:0.25] {(x-1)^2};
    \addplot[second valid,domain=0.25:1.25] {x^2+1/(8*x)};

  \end{groupplot}
  \coordinate (phiTop) at (current bounding box.north);
  \coordinate (phiBottom) at (current bounding box.south);
\end{tikzpicture}}
  \caption{Typical worst-case upper bound functions $\Phi_k(\alpha,E_k)=h_k^2\alpha^2+
  \max_{0\le e\le E_k}\{e-2\alpha c e^{p/2}\}$ with $c=E_k=h_k=1$ for $1\le p\le2$ (Left) and $p=4$ (Right). Solid curves show the valid smooth branches; dashed curves extend the same formulas beyond their valid intervals. Both functions are strongly convex.}
  \label{fig:phi-strong-convexity}
\end{figure}
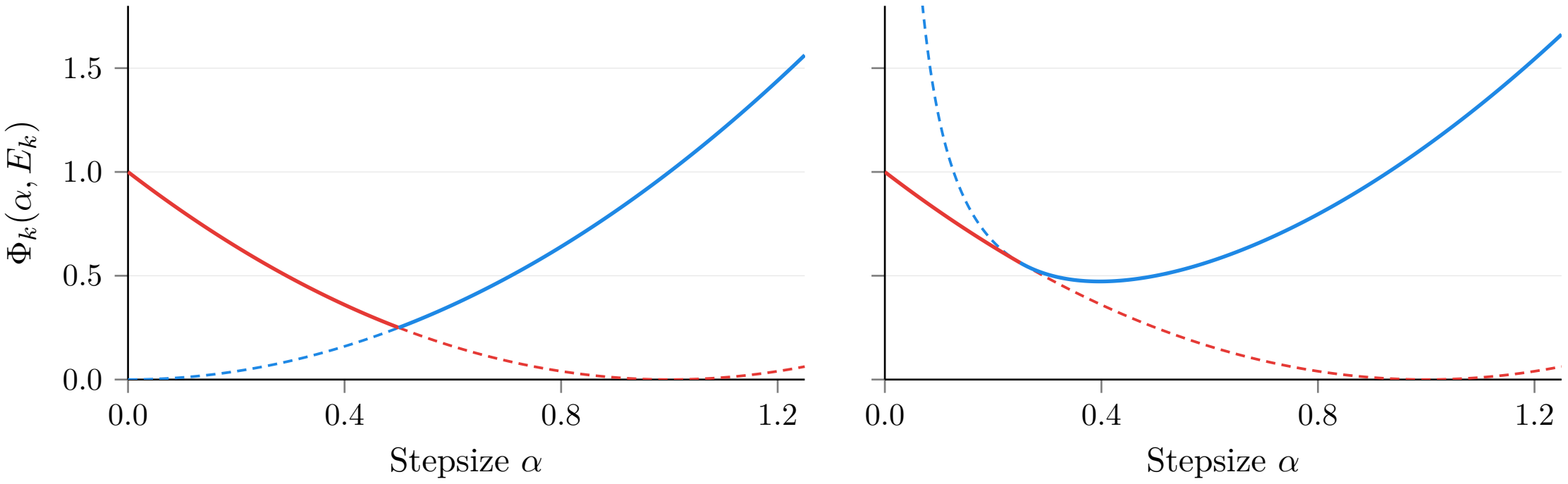

We will show in \Cref{subsection:proofs-stepsize} that the worst-case function $\Phi_k(\alpha,E_k)$ is further piecewise smooth in $\alpha$. \Cref{fig:phi-strong-convexity} illustrates its typical shapes for $1\le p\le2$ and $p=4$ that are both strongly convex. 

Before deriving the analytical solution, we provide some intuition for the minimax problem \cref{eq:stepsize-optimization}. When $h_k$ is sufficiently large, the optimal stepsize will be small enough that the inner maximization attains its worst-case value at $e=E_k$. In this case, we have
$
\Phi_k(\alpha,E_k) = \alpha^2h_k^2+ E_k-2\alpha cE_k^{p/2},
$ 
which is quadratic in $\alpha$. Minimizing this expression~gives
\begin{equation} \label{eq:distance-step-size-common}
    \alpha_{\mathrm{P}} = \frac{c E_k^{p/2}}{h_k^2}, 
\end{equation}
which has the same form as the exact distance-based Polyak-type stepsize \cref{eq:distance-based_Polyak-type_stepsize}, with the unknown squared-distance error $e_k$ replaced by its certified upper bound $E_k$. To determine when this stepsize \cref{eq:distance-step-size-common} is indeed optimal, it is convenient to introduce the following normalized quantity
\begin{equation}\label{eq:normalized-subgradient-ratio}
    \lambda_k   := \frac{h_k^2}{c^2 E_k^{p-1}},
\end{equation}
which measures the size of the observed subgradient relative to the growth parameters and the current distance certificate. For $p>1$, the certificate update \cref{eq:certificate-update} immediately implies $\lambda_k\geq 1$. For $p=1$, the same inequality follows from \Cref{lem:subgradient-distance-certificate} whenever $x_k$ is nonoptimal.

The threshold on $\lambda_k$ that determines whether $\alpha_P$ in \cref{eq:distance-step-size-common} is optimal to \cref{eq:stepsize-optimization} depends on the curvature of the function
$
    e\mapsto e-2\alpha c e^{p/2},
    \, e\in[0,E_k].
$ 
Since this curvature changes at $p=2$, we consider the regimes $1\leq p\leq 2$ and $p>2$ separately. For $1\leq p\leq 2$, the function is convex in $e$, so its maximum over $[0,E_k]$ is attained at one of the endpoints. This characterization yields the following explicit stepsize and
certificate update.
\begin{subequations}
\begin{proposition}[Explicit \DCSG{} update for $1\leq p\leq 2$]
\label{thm:dc-sg-certificate-1-leq-p-leq-2}
    Suppose that $1\leq p\leq 2$, and consider a nonoptimal iterate $x_k$ with $E_k>0$. Then the distance-certified stepsize in \cref{eq:stepsize-optimization} is unique, and the
    corresponding update has the following two cases.
    \begin{enumerate}
        \item If $\lambda_k>2$ then
        \begin{equation}\label{eq:dc-sg-update-large-lambda}
            \alpha_{\texttt{cert},k}
            =
            \frac{c E_k^{p/2}}{h_k^2},
            \qquad
            U_{k+1}
            =
            E_k-\frac{c^2E_k^p}{h_k^2}
            =
            \left(1-\frac{1}{\lambda_k}\right)E_k.
        \end{equation}

        \item If $1\leq\lambda_k\leq 2$, then
        \begin{equation}\label{eq:dc-sg-update-small-lambda}
            \alpha_{\texttt{cert},k}
            =
            \frac{E_k^{1-p/2}}{2c},
            \qquad
            U_{k+1}
            =
            \frac{h_k^2E_k^{2-p}}{4c^2}
            =
            \frac{\lambda_k}{4}E_k.
        \end{equation}
    \end{enumerate}
\end{proposition}
\end{subequations}

When $p> 2$, the function $ e\mapsto e-2\alpha c e^{p/2}, \, e\in[0,E_k].$  is concave in $e$, and the inner maximization may attain its maximum either at the endpoint $e=E_k$ or at an interior point. This leads to the following explicit distance-certified updates.

\begin{subequations}
\begin{proposition}[Explicit \DCSG{} update for $p> 2$]
\label{thm:dc-sg-certificate-p-geq-2}
    Suppose that $p> 2$, and consider a nonoptimal iterate $x_k$ with $E_k>0$. Then the distance-certified stepsize in \cref{eq:stepsize-optimization} is unique, and the corresponding update has the following two cases.
    \begin{enumerate}
        \item If $\lambda_k>p$, then
        \begin{equation}\label{eq:dc-sg-update-p-geq-2-large-lambda}
            \alpha_{\texttt{cert},k} = \frac{cE_k^{p/2}}{h_k^2},
            \qquad
            U_{k+1} = E_k-\frac{c^2E_k^p}{h_k^2}=\left(1-\frac{1}{\lambda_k}\right)E_k.
        \end{equation}

        \item If $1\leq\lambda_k\leq p$ then
        \begin{equation}\label{eq:dc-sg-stepsize-p-geq-2-small-lambda}
            \alpha_{\texttt{cert},k} = \frac{c}{h_k^2}
            \left(
                \frac{h_k^2}{pc^2}
            \right)^{\frac{p}{2(p-1)}}
            \qquad  
            U_{k+1} =  \frac{p-1}{p}
            \left(
                \frac{h_k^2}{pc^2}
            \right)^{\frac{1}{p-1}} = \frac{p-1}{p}
            \left(
                \frac{\lambda_k}{p}
            \right)^{\frac{1}{p-1}} E_k.
        \end{equation}
    \end{enumerate}
\end{proposition}
\end{subequations}

\Cref{thm:dc-sg-certificate-p-geq-2,thm:dc-sg-certificate-1-leq-p-leq-2} show that, when $\lambda_k$ in \cref{eq:normalized-subgradient-ratio} is sufficiently large, both \cref{eq:dc-sg-update-p-geq-2-large-lambda,eq:dc-sg-update-large-lambda} reduce to the common Polyak-type stepsize \cref{eq:distance-step-size-common}. When $\lambda_k$ is smaller, however, the optimal stepsize depends on the growth exponent $p$, as
shown in \cref{eq:dc-sg-stepsize-p-geq-2-small-lambda,eq:dc-sg-update-small-lambda}. 
The proofs of \Cref{thm:dc-sg-certificate-p-geq-2,thm:dc-sg-certificate-1-leq-p-leq-2} rely on careful direct calculations, and we present them in \Cref{subsection:proofs-stepsize}.

For clarity, \Cref{tab:dc-sg-explicit-updates} summarizes the explicit stepsizes and the corresponding one-step contraction factors obtained in \Cref{thm:dc-sg-certificate-p-geq-2,thm:dc-sg-certificate-1-leq-p-leq-2}. Note that these results greatly refine the monotonicity established in \Cref{prop:certificate-monotonicity} by providing explicit contraction estimates for the distance certificates, which play an essential role in the subsequent convergence analysis in \Cref{sec:dc-sg-shrinking-rates}.

\begin{table}[t]
    \centering
    \caption{Explicit distance-certified stepsizes and certificate updates. Here,
    $\lambda_k=h_k^2/(c^2E_k^{p-1})$. The two sets of updates coincide when $p = 2$.}
    \label{tab:dc-sg-explicit-updates}
    \setlength{\belowcaptionskip}{0pt}
    \renewcommand{\arraystretch}{1.2}
    \begin{tabular}{cccc}
        \toprule
        Growth exponent
        & Regime
        & Stepsize $\alpha_{\texttt{cert},k}$
        & Certificate ratio $U_{k+1}/E_k$ \\
        \midrule
        $1\leq p\leq 2$
        & $\lambda_k>2$
        & $\displaystyle {cE_k^{p/2}}/{h_k^2}$
        & $\displaystyle 1-{1}/{\lambda_k}$ \\

        $1\leq p\leq 2$
        & $1\leq\lambda_k\leq 2$
        & $\displaystyle {E_k^{1-p/2}}/{2c}$
        & $\displaystyle {\lambda_k}/{4}$ \\

        \midrule

        $p>2$
        & $\lambda_k>p$
        & $\displaystyle {cE_k^{p/2}}/{h_k^2}$
        & $\displaystyle 1-{1}/{\lambda_k}$ \\

        $p>2$
        & $1\leq\lambda_k\leq p$
        & $\displaystyle
            \frac{c}{h_k^2}
            \left(
                \frac{h_k^2}{pc^2}
            \right)^{\frac{p}{2(p-1)}}$
        & $\displaystyle
            \frac{p-1}{p}
            \left(
                \frac{\lambda_k}{p}
            \right)^{\frac{1}{p-1}}$ \\

        \bottomrule
    \end{tabular}
\end{table}

\subsection{Proofs of \Cref{thm:dc-sg-certificate-p-geq-2,thm:dc-sg-certificate-1-leq-p-leq-2} and stepsize interpretation} \label{subsection:proofs-stepsize}

We here prove the explicit formulas in \Cref{thm:dc-sg-certificate-p-geq-2,thm:dc-sg-certificate-1-leq-p-leq-2} and then discuss their interpretation.

\begin{proof}[Proof of
\Cref{thm:dc-sg-certificate-1-leq-p-leq-2}]
    Suppose that $1\leq p\leq 2$. Since  $
        e\mapsto e-2\alpha c e^{p/2}
    $ is convex on $[0,E_k]$, its maximum is attained at an endpoint. Therefore, the inner maximization \cref{eq:maximum-function} becomes 
    $
        M_k(\alpha) = \max\left\{ 0,\, E_k-2\alpha cE_k^{p/2} \right\}.
    $ 
    Define  the breakpoint stepsize
$
    \alpha_{\mathrm{bp},k} := \frac{E_k^{1-p/2}}{2c},
$ 
at which $E_k-2\alpha_{\mathrm{bp},k}cE_k^{p/2}=0.$ 
    It follows that the worst-case one-step upper bound has two branches 
    \[
        \Phi_k(\alpha,E_k)
        =
        \begin{cases}
            \alpha^2h_k^2
            +
            E_k
            -
            2\alpha cE_k^{p/2},
            & 0\leq\alpha\leq\alpha_{\mathrm{bp},k},\\[1ex]
            \alpha^2h_k^2,
            & \alpha\geq\alpha_{\mathrm{bp},k}.
        \end{cases}
    \]

    The unconstrained minimizer of the first branch is
    $
        \widehat{\alpha}_k := \frac{cE_k^{p/2}}{h_k^2}.
    $ As expected, this coincides with the Polyak-type stepsize \cref{eq:distance-step-size-common}. The second branch is strictly increasing. Consequently, the unique minimizer is
    $
        \alpha_{\texttt{cert},k} = \min\left\{ \alpha_{\mathrm{bp},k},\widehat{\alpha}_k\right\}.
    $ 
    It is easy to verify that  $\widehat{\alpha}_k<\alpha_{\mathrm{bp},k} \,\Longleftrightarrow\, h_k^2>2c^2E_k^{p-1} \,\Longleftrightarrow\, \lambda_k>2. $  Thus, when $\lambda_k>2$, the unconstrained minimizer $\widehat{\alpha}_k$ lies in $[0,\alpha_{\mathrm{bp},k}]$ and is therefore the optimal stepsize. When $1\leq\lambda_k\leq2$, the minimum is attained at the breakpoint $\alpha_{\mathrm{bp},k}$. 
    Substituting the corresponding stepsizes into $\Phi_k(\alpha,E_k)$ gives
    \cref{eq:dc-sg-update-large-lambda,eq:dc-sg-update-small-lambda}.
\end{proof}

\begin{proof}[Proof of
\Cref{thm:dc-sg-certificate-p-geq-2}]
    Consider $p>2$. In this case, the function
    $
        e\mapsto e-2\alpha c e^{p/2}
    $ 
    is concave on $[0,E_k]$. Its derivative with respect to $e$ is $  1-\alpha cp e^{(p-2)/2}.$ 
    Define the breakpoint stepsize
    $
        \alpha_{\mathrm{bp},k}
        :=
        \frac{E_k^{1-p/2}}{pc} >0.
    $ 
    If $0\leq\alpha\leq\alpha_{\mathrm{bp},k}$, the inner maximum is attained at $e=E_k$. If $\alpha\geq\alpha_{\mathrm{bp},k}$, an interior maximizer is given by
    $
        e_\alpha := \left(\frac{1}{\alpha cp}\right)^{\frac{2}{p-2}}.
    $ 
    Using the identity
    $
        \alpha cp e_\alpha^{(p-2)/2}=1,
    $
    we obtain
    \[
        M_k(\alpha)
        =
        \begin{cases}
            E_k-2\alpha cE_k^{p/2},
            & 0\leq\alpha\leq\alpha_{\mathrm{bp},k},\\[1ex]
            \displaystyle
            \frac{p-2}{p}
            \left(\frac{1}{\alpha cp}\right)^{\frac{2}{p-2}},
            & \alpha\geq\alpha_{\mathrm{bp},k}.
        \end{cases}
    \]
    Therefore, the worst-case one-step upper bound also has two branches:
    \[
        \Phi_k(\alpha,E_k)
        =
        \begin{cases}
            \alpha^2h_k^2+E_k-2\alpha cE_k^{p/2},
            & 0\leq\alpha\leq\alpha_{\mathrm{bp},k},\\[1ex]
            \displaystyle
            \alpha^2h_k^2 + \frac{p-2}{p}\left(\frac{1}{\alpha cp}\right)^{\frac{2}{p-2}},
            & \alpha\geq\alpha_{\mathrm{bp},k}.
        \end{cases}
    \]

    The unconstrained minimizer of the first branch is
    $
        \widehat{\alpha}_k  := \frac{cE_k^{p/2}}{h_k^2},
    $ 
    which again coincides with the Polyak-type stepsize
    \cref{eq:distance-step-size-common}. Moreover, we have 
    $
        \widehat{\alpha}_k<\alpha_{\mathrm{bp},k}
        \,\Longleftrightarrow\,
        h_k^2>pc^2E_k^{p-1}
        \,\Longleftrightarrow\,
        \lambda_k>p.
    $ 
    Thus, when $\lambda_k>p$, the unconstrained minimizer
    $\widehat{\alpha}_k$ lies in $[0,\alpha_{\mathrm{bp},k}]$ and is thus the
    optimal stepsize.

    The unconstrained minimizer of the second branch is
    $
        \widetilde{\alpha}_k
        := \frac{c}{h_k^2}
            \left(
                \frac{h_k^2}{pc^2}
            \right)^{\frac{p}{2(p-1)}}.
    $ 
    A direct calculation gives
    $
        \frac{\widetilde{\alpha}_k}
             {\alpha_{\mathrm{bp},k}}
        =
        \left(\frac{p}{\lambda_k}\right)^{
            \frac{p-2}{2(p-1)}
        }.
    $ 
    Hence, when $1\leq\lambda_k\leq p$, we have $\widetilde{\alpha}_k\geq\alpha_{\mathrm{bp},k}$, so $\widetilde{\alpha}_k$ lies in $[\alpha_{\mathrm{bp},k},\infty)$ and is the optimal stepsize. At $\lambda_k=p$, both branch minimizers coincide with $\alpha_{\mathrm{bp},k}$. Substituting the corresponding stepsizes into $\Phi_k(\alpha,E_k)$ yields \cref{eq:dc-sg-stepsize-p-geq-2-small-lambda,eq:dc-sg-update-p-geq-2-large-lambda}
\end{proof}

\begin{figure}[t]
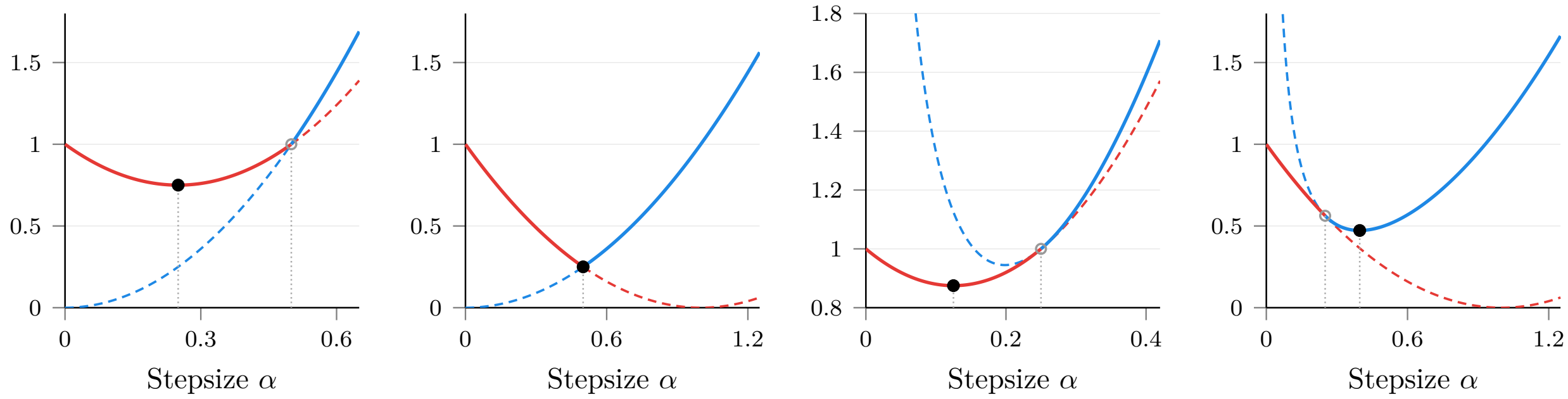

    \centering
    \begingroup
\captionsetup[subfigure]{font=small,justification=centering,
    singlelinecheck=false,skip=5pt}
\begin{subfigure}[t]{0.24\linewidth}
    \centering
    \resizebox{\linewidth}{!}{\begin{tikzpicture}
  \input{figures/phi_table2_style.tex}
  \path[use as bounding box] (-0.80,-0.90) rectangle (3.65,3.58);
  \begin{axis}[dc table two,
      xmin=0,xmax=0.65,ymin=0,ymax=1.8,
      xtick={0,0.3,0.6},ytick={0,0.5,1,1.5}
  ]
    \addplot[endpoint extension,domain=0.5:0.65] {4*x^2+1-2*x};
    \addplot[second extension,domain=0:0.5] {4*x^2};
    \addplot[endpoint valid,domain=0:0.5] {4*x^2+1-2*x};
    \addplot[second valid,domain=0.5:0.65] {4*x^2};
    \addplot[guide] coordinates {(0.25,0) (0.25,0.75)};
    \addplot[guide] coordinates {(0.5,0) (0.5,1)};
    \addplot[breakpoint] coordinates {(0.5,1)};
    \addplot[selected] coordinates {(0.25,0.75)};
  \end{axis}
\end{tikzpicture}}
    \caption{Polyak-type\\$1\le p\le2,\;\lambda_k=4>2$\\$\alpha_{\mathrm{cert},k}=\widehat\alpha_k=1/4$}
    \label{fig:phi-regimes-a}
\end{subfigure}\hfill%
\begin{subfigure}[t]{0.24\linewidth}
    \centering
    \resizebox{\linewidth}{!}{\begin{tikzpicture}
  \input{figures/phi_table2_style.tex}
  \path[use as bounding box] (-0.80,-0.90) rectangle (3.65,3.58);
  \begin{axis}[dc table two,
      xmin=0,xmax=1.25,ymin=0,ymax=1.8,
      xtick={0,0.6,1.2},ytick={0,0.5,1,1.5}
  ]
    \addplot[endpoint extension,domain=0.5:1.25] {x^2+1-2*x};
    \addplot[second extension,domain=0:0.5] {x^2};
    \addplot[endpoint valid,domain=0:0.5] {x^2+1-2*x};
    \addplot[second valid,domain=0.5:1.25] {x^2};
    \addplot[guide] coordinates {(0.5,0) (0.5,0.25)};
    \addplot[breakpoint] coordinates {(0.5,0.25)};
    \addplot[selected] coordinates {(0.5,0.25)};
  \end{axis}
\end{tikzpicture}}
    \caption{Capped\\$1\le p\le2,\;\lambda_k=1\le2$\\$\alpha_{\mathrm{cert},k}=\alpha_{\mathrm{bp},k}=1/2$}
    \label{fig:phi-regimes-b}
\end{subfigure}\hfill%
\begin{subfigure}[t]{0.24\linewidth}
    \centering
    \resizebox{\linewidth}{!}{\begin{tikzpicture}
  \input{figures/phi_table2_style.tex}
  \path[use as bounding box] (-0.80,-0.90) rectangle (3.65,3.58);
  \begin{axis}[dc table two,
      xmin=0,xmax=0.42,ymin=0.8,ymax=1.8,
      xtick={0,0.2,0.4},ytick={0.8,1,1.2,1.4,1.6,1.8}
  ]
    \addplot[endpoint extension,domain=0.25:0.42] {8*x^2+1-2*x};
    \addplot[second extension,domain=0.02:0.25] {8*x^2+1/(8*x)};
    \addplot[endpoint valid,domain=0:0.25] {8*x^2+1-2*x};
    \addplot[second valid,domain=0.25:0.42] {8*x^2+1/(8*x)};
    \addplot[guide] coordinates {(0.125,0.8) (0.125,0.875)};
    \addplot[guide] coordinates {(0.25,0.8) (0.25,1)};
    \addplot[breakpoint] coordinates {(0.25,1)};
    \addplot[selected] coordinates {(0.125,0.875)};
  \end{axis}
\end{tikzpicture}}
    \caption{Polyak-type\\$p=4,\;\lambda_k=8>p$\\$\alpha_{\mathrm{cert},k}=\widehat\alpha_k=1/8$}
    \label{fig:phi-regimes-c}
\end{subfigure}\hfill%
\begin{subfigure}[t]{0.24\linewidth}
    \centering
    \resizebox{\linewidth}{!}{\begin{tikzpicture}
  \input{figures/phi_table2_style.tex}
  \pgfmathsetmacro{\dcInteriorStep}{pow(1/16,1/3)}
  \pgfmathsetmacro{\dcInteriorValue}{\dcInteriorStep^2+1/(8*\dcInteriorStep)}
  \path[use as bounding box] (-0.80,-0.90) rectangle (3.65,3.58);
  \begin{axis}[dc table two,
      xmin=0,xmax=1.25,ymin=0,ymax=1.8,
      xtick={0,0.6,1.2},ytick={0,0.5,1,1.5}
  ]
    \addplot[endpoint extension,domain=0.25:1.25] {x^2+1-2*x};
    \addplot[second extension,domain=0.02:0.25] {x^2+1/(8*x)};
    \addplot[endpoint valid,domain=0:0.25] {x^2+1-2*x};
    \addplot[second valid,domain=0.25:1.25] {x^2+1/(8*x)};
    \addplot[guide] coordinates {(0.25,0) (0.25,0.5625)};
    \addplot[guide] coordinates {(\dcInteriorStep,0) (\dcInteriorStep,\dcInteriorValue)};
    \addplot[breakpoint] coordinates {(0.25,0.5625)};
    \addplot[selected] coordinates {(\dcInteriorStep,\dcInteriorValue)};
  \end{axis}
\end{tikzpicture}}
    \caption{Interior-maximizer\\$p=4,\;\lambda_k=1\le p$\\$\alpha_{\mathrm{cert},k}=\widetilde\alpha_k\approx0.397$}
    \label{fig:phi-regimes-d}
\end{subfigure}
\endgroup
    \caption{The four minimizer regimes in \Cref{tab:dc-sg-branch-summary}. Each panel plots $\Phi_k(\alpha,E_k)$ against the stepsize $\alpha$ where $c=E_k=1$ and $h_k^2=\lambda_k$.  Solid curves show the valid branches; dashed curves extend their formulas outside the valid intervals. Filled dots mark the certified minimizers $\alpha_{\mathrm{cert},k}$; open circles mark the breakpoints $\alpha_{\mathrm{bp},k}$.  The symbols $\widehat\alpha_k$ and $\widetilde\alpha_k$ denote the minimizers of the first and second smooth branches, respectively.}
    \label{fig:phi-regimes}
\end{figure}

The explicit updates in \cref{thm:dc-sg-certificate-1-leq-p-leq-2,thm:dc-sg-certificate-p-geq-2} share a common two-branch structure, which is summarized in \Cref{tab:dc-sg-branch-summary}. 
In particular, for all $p\geq1$, the first branch is the Polyak-type update  \cref{eq:dc-sg-update-large-lambda,eq:dc-sg-update-p-geq-2-large-lambda}. This stepsize is selected when $\lambda_k$ is sufficiently large, namely, when $\lambda_k>2$ for $1\leq p\leq2$ and when $\lambda_k>p$ for $p>2$.  
The complementary branch depends on the curvature of $e-2\alpha c e^{p/2}$. For $1\leq p\leq2$, this function is convex in $e$, and its maximum over $[0,E_k]$ is attained at an endpoint. When $1\leq\lambda_k\leq2$, the optimal stepsize is capped at the breakpoint $\alpha_{\mathrm{bp},k} := \frac{E_k^{1-p/2}}{2c}$, leading to \cref{eq:dc-sg-update-small-lambda}. 
We refer to this update as the \textit{capped branch}. For $p>2$, the function $e-2\alpha c e^{p/2}$ is concave, and the inner maximum may be attained at an interior stationary point. When $1\leq\lambda_k\leq p$, this gives \cref{eq:dc-sg-stepsize-p-geq-2-small-lambda}. We refer to this update as the \textit{interior-maximizer branch}. \Cref{fig:phi-regimes} illustrates these four different  minimizer regimes in \Cref{tab:dc-sg-branch-summary}. 

\begin{table}[t]
    \centering
    \caption{Summary of the two branches in the explicit \DCSG{} updates; $\lambda_k$ is defined in \Cref{eq:normalized-subgradient-ratio}.}
    \label{tab:dc-sg-branch-summary}
    \renewcommand{\arraystretch}{1.2}
    \begin{tabular}{
        @{}c
        >{\raggedright\arraybackslash}p{0.34\linewidth}
        >{\raggedright\arraybackslash}p{0.39\linewidth}
        @{}
    }
        \toprule
        \textbf{Range of $p$}
        &
        \textbf{Common branch}
        \newline
        $\lambda_k>\max\{2,p\}$
        &
        \textbf{Complementary branch}
        \newline
        $1\leq\lambda_k\leq\max\{2,p\}$
        \\
        \midrule

        $1\leq p\leq2$ & Polyak-type branch; see
        \cref{eq:dc-sg-update-large-lambda}. &
        Capped branch; see \cref{eq:dc-sg-update-small-lambda}.
        \\
        $p>2$ & Polyak-type branch; see~\cref{eq:dc-sg-update-p-geq-2-large-lambda}. & Interior-maximizer branch; see~\cref{eq:dc-sg-stepsize-p-geq-2-small-lambda}.
        \\
        \bottomrule
    \end{tabular}
    \vspace{8pt}
\end{table}

Despite their different forms, the distance certificate is always contracted. The explicit contraction factors derived in \cref{thm:dc-sg-certificate-1-leq-p-leq-2,thm:dc-sg-certificate-p-geq-2} form the basis of the convergence analysis in \Cref{sec:dc-sg-shrinking-rates}.

\section{Convergence rates of \DCSG} \label{sec:dc-sg-shrinking-rates}

In this section, we establish the convergence rates of our proposed \DCSG{} method in \Cref{alg:dc-sg}. We also compare our \DCSG{} method with the state-of-the-art results. 

\subsection{Convergence guarantees under global $L$-Lipschitz continuity}
Suppose that $f$ is globally $L$-Lipschitz continuous. Then every subgradient satisfies $
    h_k=\|g_k\|\leq L.$ For convenience, define
\begin{equation}\label{eq:kappa-q-definition}
    \kappa:=\frac{L}{c},  \qquad q:=\max\{2,p\}.
\end{equation}
Based on \Cref{thm:dc-sg-certificate-1-leq-p-leq-2,thm:dc-sg-certificate-p-geq-2}, the following result provides a unified one-step decrease for all~$p\geq1$.

\begin{lemma}[One-step certificate decrease]
\label{thm:one-step-certificate-lipschitz}
    Let $f$ be convex and globally $L$-Lipschitz continuous. Suppose that $f$ satisfies the $(p,c)$-H{\"o}lder growth \cref{eq:holderian-error-bound}. For each iteration of \DCSG, we~have
    \begin{equation}\label{eq:one-step-certificate-lipschitz}
        E_{k+1} \leq   E_k    -  \min\left\{
            \frac{E_k}{q},
            \frac{E_k^p}{\kappa^2}
        \right\}.
    \end{equation}
\end{lemma}

\begin{proof}
    Suppose that $\lambda_k > \max\{2,p\}$, i.e., the Polyak branch is active. By
    \Cref{thm:dc-sg-certificate-1-leq-p-leq-2,thm:dc-sg-certificate-p-geq-2}, its propagated certificate satisfies
    $
        U_{k+1} = E_k-\frac{c^2E_k^p}{h_k^2}.
    $ 
    Using $h_k\leq L$ gives
    \begin{equation} \label{eq:polyak-update}
        U_{k+1}  \leq  E_k-\frac{c^2E_k^p}{L^2} = E_k-\frac{E_k^p}{\kappa^2}.
    \end{equation}

    Suppose that $\lambda_k \leq \max\{2,p\}$, i.e., the complementary branch is active. When $1\leq p\leq2$, the capped update \cref{eq:dc-sg-update-small-lambda} satisfies
    $
        U_{k+1}\leq{E_k}/{2}.
    $ 
    When $p>2$, the interior-maximizer update \cref{eq:dc-sg-stepsize-p-geq-2-small-lambda} satisfies
    $
        U_{k+1}\leq
        \left(1-\frac{1}{p}\right)E_k.
    $ 
    Since $q=\max\{2,p\}$, both cases imply
     \begin{equation} \label{eq:complementary-update}
        U_{k+1} \leq E_k-\frac{E_k}{q}.
    \end{equation}
    Therefore, regardless of the value of $\lambda_k$, combining \cref{eq:polyak-update} with \cref{eq:complementary-update} leads to  
    \[
        U_{k+1} \leq  E_k - \min\left\{
            \frac{E_k}{q},
            \frac{E_k^p}{\kappa^2}
        \right\}.
    \]
    By construction \cref{eq:certificate-update}, the certificate update guarantees $E_{k+1}\leq U_{k+1}$,
    the bound \cref{eq:one-step-certificate-lipschitz} follows.
\end{proof}

The recursion \cref{eq:one-step-certificate-lipschitz} directly yields global convergence of the certificate sequence $\{E_k\}$. We first consider the case $p>1$. Throughout this section, define the threshold
\begin{equation}\label{eq:lipschitz-certificate-threshold}
    E_{\mathrm{th}}
    :=
    \left(\frac{\kappa^2}{q}\right)^{\frac{1}{p-1}},
\end{equation}
at which we have $ {E_k}/{q}={E_k^p}/{\kappa^2}$. Above this threshold, the certificate decreases geometrically until it enters the local region; within the local region, it converges at a
sublinear rate.

\begin{theorem}[Sublinear convergence when $p>1$]
\label{thm:polynomial-certificate-decay-lipschitz}
    Let $f$ be convex and globally $L$-Lipschitz continuous. Suppose that $f$ satisfies the $(p,c)$-H{\"o}lder growth \cref{eq:holderian-error-bound} with $p>1$. Define $E_{\mathrm{th}}$ in \cref{eq:lipschitz-certificate-threshold}.
    \vspace{-1mm}
    \begin{subequations}
    \begin{itemize}
        \item  If $E_0\leq E_{\mathrm{th}}$, then, for every $k\geq0$,
    \begin{equation}\label{eq:local-polynomial-rate-lipschitz}
        E_k
        \leq
        \left(
            E_0^{1-p}
            +
            \frac{(p-1)k}{\kappa^2}
        \right)^{-\frac{1}{p-1}}.
    \end{equation}

    \item If $E_0>E_{\mathrm{th}}$, define
    $   
        K_{0}
        :=
        \left\lceil
            \frac{\log(E_0/E_{\mathrm{th}})}
                 {\log(q/(q-1))}
        \right\rceil.
    $ 
    Then $E_{K_{0}}\leq E_{\mathrm{th}}$, and, for every
    $k\geq K_{0}$,
    \begin{equation}\label{eq:global-polynomial-rate-lipschitz}
        E_k \leq \left(
            \frac{
                q+(p-1)(k-K_{0})
            }{\kappa^2}
        \right)^{-\frac{1}{p-1}}.
    \end{equation}
    \end{itemize}
    \end{subequations}
\end{theorem}

\begin{proof}
    Suppose first that $E_0\leq E_{\mathrm{th}}$. Since the certificate
    sequence is nonincreasing, $E_k\leq E_{\mathrm{th}}$ for all
    $k\geq0$. Hence, we have
    $
        \frac{E_k^p}{\kappa^2}
        \leq \frac{E_k}{q},
    $
    and \cref{eq:one-step-certificate-lipschitz} reduces to
    \[
        E_{k+1}
        \leq
        E_k\left(
            1-\frac{E_k^{p-1}}{\kappa^2}
        \right).
    \]
    If $E_k=0$ for some $k$, \cref{eq:local-polynomial-rate-lipschitz} follows immediately from monotonicity. Otherwise, applying Bernoulli's inequality gives
    \[
    \begin{aligned}
        E_{k+1}^{1-p}
        \geq
        E_k^{1-p}
        \left(
            1-\frac{E_k^{p-1}}{\kappa^2}
        \right)^{1-p} \geq
        E_k^{1-p}+\frac{p-1}{\kappa^2}.
    \end{aligned}
    \]
    Iterating this inequality from $0$ to $k-1$ yields \cref{eq:local-polynomial-rate-lipschitz}.

    Suppose next that $E_0>E_{\mathrm{th}}$. Then we have
    $
        {E_k^p}/{\kappa^2}
        >
        {E_k}/{q},
    $ 
    and hence \cref{eq:one-step-certificate-lipschitz} reduces to
    $
        E_{k+1}
        \leq
        \left(1-\frac{1}{q}\right)E_k.
    $ 
    Therefore, the certificate either enters the local region before
    $K_0$ or satisfies
    \[
        E_{K_0}
        \leq
        E_0\left(1-\frac{1}{q}\right)^{K_0}
        \leq E_{\mathrm{th}}.
    \]
    Applying \cref{eq:local-polynomial-rate-lipschitz} from iteration
    $K_0$ and using
    $
        E_{K_0}^{1-p}  \geq  E_{\mathrm{th}}^{1-p} = {q}/{\kappa^2}
    $ 
    proves \cref{eq:global-polynomial-rate-lipschitz}.
\end{proof}

\Cref{thm:polynomial-certificate-decay-lipschitz} immediately yields the following iteration complexity.

\begin{corollary}[Iteration complexity for $p>1$]
\label{cor:iteration-complexity-lipschitz}
    Under the setting of
    \Cref{thm:polynomial-certificate-decay-lipschitz}, let $0<\epsilon<E_0$. \DCSG{} in \Cref{alg:dc-sg}  obtains an iterate satisfying  $\dist(x_k,X^\star)^2\leq\epsilon$ within 
    $
        O\left(\kappa^2\epsilon^{-(p-1)}\right)
    $ 
    iterations. More precisely, the following statements hold. 
\begin{itemize}
    \item if $E_0\leq E_{\mathrm{th}}$, then $E_k\leq\epsilon$ whenever
    $
        k \geq \left\lceil
            \frac{\kappa^2}{p-1}
            \left(
                \epsilon^{1-p}-E_0^{1-p}
            \right)
        \right\rceil.
   $

   \item  if $E_0>E_{\mathrm{th}}$ and $0<\epsilon<E_{\mathrm{th}}$, then
    $E_k\leq\epsilon$ whenever
    $ k \geq K_0  +\left\lceil
            \frac{\kappa^2\epsilon^{1-p}-q}{p-1}
        \right\rceil.
   $ 
\end{itemize}
\end{corollary}

When $p=1$, \cref{eq:one-step-certificate-lipschitz} immediately gives linear convergence. We present the guarantee below. 

\begin{theorem}[Linear convergence for $p=1$]
\label{thm:geometric-certificate-decay-p-equals-1}
    Let $f$ be convex and globally $L$-Lipschitz continuous. Suppose that $f$ satisfies the $(p,c)$-H{\"o}lder growth \cref{eq:holderian-error-bound} with $p=1$. Define
    $ \delta  := \min\left\{
            \frac{1}{2},
            \frac{1}{\kappa^2}
        \right\}.
    $ 
    Then, for every $k\geq0$,
    \begin{equation}\label{eq:geometric-rate-p-equals-1}
        E_k\leq(1-\delta)^kE_0.
    \end{equation}
    Consequently, for every $0<\epsilon<E_0$, we have
    $E_k\leq\epsilon$ whenever
    \begin{equation}\label{eq:geometric-complexity-p-equals-1}
        k
        \geq
        \left\lceil
            \frac{\log(E_0/\epsilon)}
                 {\log(1/(1-\delta))}
        \right\rceil.
    \end{equation}
\end{theorem}

\begin{proof}
    When $p=1$, we have $q=2$ by definition \cref{eq:kappa-q-definition}. Thus, \cref{eq:one-step-certificate-lipschitz} becomes
    \[
        E_{k+1}  \leq E_k-\min\left\{
            \frac{E_k}{2},
            \frac{E_k}{\kappa^2}
        \right\} = (1-\delta)E_k.
    \]
    Iterating this inequality proves \cref{eq:geometric-rate-p-equals-1}. The complexity bound follows
    by solving $(1-\delta)^kE_0\leq\epsilon$.
\end{proof}

The case $p=1$ corresponds to the classical \textit{sharp-growth} condition, also known as weak sharp minima when the solution set is not necessarily a singleton \cite{BurkeFerris1993}. It has long been recognized that
this growth structure allows suitably designed subgradient methods to achieve linear convergence; see, for example, \cite{Polyak1969,goffin1977convergence,JohnstoneMoulin2020}. Thus, \Cref{thm:geometric-certificate-decay-p-equals-1} recovers the expected linear rate using our \DCSG{} method. Note that \DCSG{} maintains a
valid current-iterate distance certificate and does not require knowledge of the optimal value $f^\star$ \cite{Polyak1969} or a predetermined stagewise stepsize schedule \cite{JohnstoneMoulin2020}. We further compare with state-of-the-art analysis in \Cref{subsection:comparison}. 

\subsection{Refined rates under shrinking subgradients}
\label{subsec:shrinking-subgradient-rates}

The results in \Cref{thm:geometric-certificate-decay-p-equals-1,thm:polynomial-certificate-decay-lipschitz} assume a global Lipschitz constant. \DCSG{} can achieve faster rates when the norms of selected subgradients decrease as the iterates approach the solution set. 

Specifically, suppose that the subgradients of $f$ satisfy the following upper bound  
\begin{equation}\label{eq:shrinking-selected-subgradients}
    \|g_k\| \leq
    L\,\dist(x_k,X^\star)^s
\end{equation}
for some $L>0$ and $s\in[0,p-1]$. When $s = 0$, it reduces to the globally Lipschitz case. Note that it is necessary to have $s\leq p-1$  since the H{\"o}lder growth condition and the subgradient inequality imply
$
    h_k \geq
    c\,\dist(x_k,X^\star)^{p-1}
$ 
at every nonoptimal iterate (see \Cref{lem:subgradient-distance-certificate}). Under this shrinking subgradient assumption \cref{eq:shrinking-selected-subgradients}, we obtain the following
refinement of the rates in \Cref{thm:geometric-certificate-decay-p-equals-1,thm:polynomial-certificate-decay-lipschitz}.

\begin{corollary}[Refined rates under shrinking subgradients] \label{cor:shrinking-subgradient-rates}
    Suppose that $f$ is convex and satisfies the $(p,c)$-H{\"o}lder growth condition and that \cref{eq:shrinking-selected-subgradients} holds. Let $\kappa=L/c$ and $q=\max\{2,p\}$. Then every \DCSG{} iteration in \Cref{alg:dc-sg} satisfies
    \begin{equation}\label{eq:one-step-shrinking-subgradients}
        E_{k+1} \leq E_k-\min\left\{ \frac{E_k}{q},\frac{E_k^{p-s}}{\kappa^2} \right\}.
    \end{equation}
Consequently, the following statements hold. 
\begin{itemize}
    \item 
    If $0\leq s<p-1$, then we have $E_k = O\left(   k^{-\frac{1}{p-1-s}} \right),$ 
    and the iteration complexity for obtaining $E_k\leq\epsilon$ is
    $ 
        O\left(
            \kappa^2\epsilon^{-(p-1-s)}
        \right).
    $

    \item    If $s=p-1$, define
    $
        \delta_s
        :=
        \min\left\{
            \frac{1}{q},
            \frac{1}{\kappa^2}
        \right\}.
    $ 
    Then we have linear convergence
    $
        E_k\leq(1-\delta_s)^kE_0.
    $
\end{itemize}
\end{corollary}

\begin{proof}
    Since $e_k\leq E_k$, \cref{eq:shrinking-selected-subgradients} gives $ h_k^2 \leq L^2E_k^s. $ Thus, on the Polyak-type branch, we have
    \[ 
    \frac{c^2E_k^p}{h_k^2} \geq \frac{c^2E_k^p}{L^2E_k^s} = \frac{E_k^{p-s}}{\kappa^2}.
    \]
    The complementary branch continues to decrease the certificate by at least $E_k/q$, similar to \cref{eq:complementary-update}. Thus, \cref{eq:one-step-certificate-lipschitz} reduces to  \cref{eq:one-step-shrinking-subgradients}.

    If $s<p-1$, applying the argument of \Cref{thm:polynomial-certificate-decay-lipschitz} with $p-1$ replaced by $p-1-s$ establish the claimed result. If $s=p-1$, then \cref{eq:one-step-shrinking-subgradients} reduces to $E_{k+1}\leq(1-\delta_s)E_k.$ 
    This finishes the proof.
\end{proof}

Whenever $s>0$, \DCSG{} achieves a strictly faster convergence rate and thus automatically adapts to shrinking subgradient norms. In particular, when $p>1$ and $s=p-1$, \DCSG{} converges linearly without requiring the sharp-growth condition. For example, the function $f(x)=\|Ax\|^p$ with $A \in \mathbb{R}^{m\times n}, p>1$ satisfies the condition \cref{eq:shrinking-selected-subgradients} with $s=p-1$, and hence \DCSG{} converges linearly on this class of problems (see our experiments in \Cref{subsection:power-objective}). A related acceleration phenomenon for adaptive Polyak-type stepsizes under higher-order growth was recently established in
\cite{davis2025gradient}.

\begin{remark}[Shrinking subgradients and norm-aware stepsizes] \label{remark:norm-aware-stepsizes}
    The refined rate in \Cref{cor:shrinking-subgradient-rates} is not specific to our \DCSG{} method.
    More generally, such an improvement can arise in subgradient methods whose one-step progress bounds use the observed norm $h_k$ and do not replace it by the global upper bound $L$; see also \Cref{remark:upper-bound}. For example, the classical Polyak stepsize \cref{eq:polyak-stepsize} satisfies
    \begin{equation} \label{eq:polyak-recursion}
    \begin{aligned}
        e_{k+1}
        \leq
        e_k-\frac{(f(x_k)-f^\star)^2}{h_k^2}
        \leq
        e_k-\frac{c^2e_k^p}{h_k^2}
        \leq
        e_k-\frac{1}{\kappa^2}e_k^{p-s},
    \end{aligned}
    \end{equation}
    where the last inequality uses \cref{eq:shrinking-selected-subgradients}. Thus, the Polyak stepsize enjoys the same refined polynomial exponent when $s<p-1$ and linear convergence when $s=p-1$, similar to \Cref{cor:shrinking-subgradient-rates}. 
    The key to obtaining \cref{eq:polyak-recursion,eq:one-step-shrinking-subgradients} is to exploit \cref{eq:shrinking-selected-subgradients} directly and not to relax it using the uniform bound $h_k\leq L$.  In contrast, standard stagewise and geometrically decaying stepsize schemes \cite{JohnstoneMoulin2020,Wang2025,goffin1977convergence,YangLin2018} use the global constant $L$ in their stepsize design. Such methods do not automatically adapt to shrinking subgradient norms and may fail to attain the refined~rates (see our experiments in \Cref{section:experiments}). \hfill $\square$
\end{remark}

\subsection{Comparison with the state of the art}
\label{subsection:comparison}

We compare \DCSG{} with several closely related subgradient methods. First, the convergence rates in \Cref{thm:polynomial-certificate-decay-lipschitz,thm:geometric-certificate-decay-p-equals-1} are faster than the classical subgradient method. Under uniformly bounded subgradients, the standard analysis with a suitably chosen diminishing stepsize of order $\alpha_k=\mathcal{O}(1/\sqrt{k})$ guarantees 
$
    f(\hat{x}_k)-f^\star = O\left(k^{-1/2}\right)
$ 
for a \textit{suitable averaged iterate} $\hat{x}_k$. Combining this bound with the H{\"o}lder growth \cref{eq:holderian-error-bound} gives
$
    \dist(\hat{x}_k,X^\star) = O\left(k^{-\frac{1}{2p}}\right). 
$ 
By contrast, \Cref{thm:polynomial-certificate-decay-lipschitz,thm:geometric-certificate-decay-p-equals-1} establish the \textit{last-iterate} guarantee
\begin{equation} \label{eq:optimal-rate}
    \dist(x_k,X^\star)\leq \sqrt{E_k} = O\left(k^{-\frac{1}{2(p-1)}}\right)
\end{equation}
when $p>1$, and linear convergence when $p=1$. Thus, \DCSG{} strictly improves the polynomial exponent for every $p>1$ and achieves a linear rate in the sharp-growth case. Moreover, these guarantees hold directly for the last iterate $x_k$ and do not require averaging.

As discussed in \Cref{subsection:stepsize-design}, our distance-certified stepsize \cref{eq:stepsize-optimization} is directly motivated by the classical Polyak stepsize \cref{eq:polyak-stepsize}. Indeed, both stepsizes aim to minimize the one-step upper bound shown in \Cref{lem:distance-recursion}, but they use different information. The Polyak stepsize requires exact knowledge of the optimal value $f^\star$, and our \DCSG{} uses the growth parameters $(p,c)$ and an initial certificate $U_0$. As shown in \cref{eq:polyak-recursion}, the Polyak stepsize attains the same asymptotic rate as \DCSG{} under H{\"o}lder growth and similarly benefits from shrinking subgradients. Recently, \cite[Theorems~1 and~2]{Wang2025} established matching upper and lower bounds showing that the Polyak stepsize is minimax optimal for decreasing the distance to the solution set for a class of globally Lipschitz convex functions satisfying a general error bound. When specialized to $(p,c)$-H{\"o}lder growth, the guarantees in \Cref{thm:polynomial-certificate-decay-lipschitz,thm:geometric-certificate-decay-p-equals-1} also attain the corresponding minimax-optimal rate asymptotically.

There are several other stepsize schedules for subgradient methods based on growth parameters and a global Lipschitz constant, including \cite{goffin1977convergence,JohnstoneMoulin2020,YangLin2018,Wang2025}.
For the sharp-growth case $p=1$, Goffin \cite{goffin1977convergence} established linear convergence using a geometrically decaying stepsize. Johnstone and Moulin \cite{JohnstoneMoulin2020} developed related decaying and descending-stairs stepsizes for general H{\"o}lder growth. In their notation, $\theta=1/p$, and their methods attain the rate in \cref{eq:optimal-rate} when $p>1$ and linear convergence when $p=1$. Yang and Lin \cite{YangLin2018} obtain comparable improvements under a local error-bound condition through a restarted subgradient method, but their complexity generally includes an additional logarithmic factor and their method requires averaging within each stage. As mentioned above, Wang \cite{Wang2025} proved that the Polyak stepsize is minimax optimal for globally Lipschitz convex functions satisfying a general error bound and, under an additional regularity condition, established the same conclusion for a predetermined decaying scheme called \texttt{DecayGD}.

All these geometric, descending-stairs, restarted, and \texttt{DecayGD} schemes \cite{goffin1977convergence,JohnstoneMoulin2020,YangLin2018,Wang2025} \textit{predetermine} their stepsizes based on a uniform subgradient bound, growth information, and an initial distance bound (see \Cref{tab:method-information} for a comparison). In contrast, \DCSG{} recomputes its stepsize adaptively at every iteration from the current certificate $E_k$ and the observed subgradient norm $h_k$. In particular, neither a global Lipschitz constant nor a predetermined stage length is required in \DCSG; the Lipschitz constant enters only its worst-case analysis.  Moreover, retaining the observed norm $h_k$ allows \DCSG{} to achieve the refined rates in \Cref{cor:shrinking-subgradient-rates} when the selected subgradients diminish near $X^\star$. The existing guarantees for the predetermined schemes above \cite{goffin1977convergence,JohnstoneMoulin2020,YangLin2018,Wang2025} are based on uniform subgradient bounds and do not capture this additional improvement; also see \Cref{remark:norm-aware-stepsizes}.

\section{Numerical experiments} \label{section:experiments}

In this section, we present two experiments to demonstrate the theoretical findings in this work. We compare our proposed \DCSG{} in \Cref{alg:dc-sg} with the classical Polyak subgradient method~\cite{Polyak1969}, the restarted subgradient method (RSG)~\cite{YangLin2018}, the descending-stairs subgradient method (DS-SG)~\cite{JohnstoneMoulin2020}, and the decaying-stepsize method
(DecayGD)~\cite{Wang2025}.

\Cref{tab:method-information} summarizes the problem information used in these stepsize schedules. Our \DCSG{} uses $p$, $c$, and $U_0$, and adapts its stepsize to the observed subgradient norm without
requiring $f^\star$ or $L$. The Polyak method requires the exact optimal value $f^\star$ and serves as a baseline. For comparison, RSG, DS-SG, and DecayGD use $p$, $c$, $L$, and $U_0$ to determine their schedules. RSG and DS-SG employ stagewise schedules, whereas DecayGD uses a predetermined decaying sequence.
We include DS-SG only when $1 \leq p \leq 2$, the range covered by its convergence guarantees.

\begin{table}[t]
    \centering
    \caption{Problem information and stepsize schedules. Here, $f^\star$ is the optimal value, $p$ and $c$ are the H\"older growth parameters, $L$ is a uniform bound on the subgradient norms, and
$U_0 \geq \operatorname{dist}(x_0,X^\star)^2$ is an initial squared-distance bound. A checkmark ($\checkmark$) indicates that the corresponding quantity is used in the method.}
    \label{tab:method-information}
    \small
    \setlength{\tabcolsep}{6pt}
    \renewcommand{\arraystretch}{1.15}
    \begin{tabular}{@{}lcccccll@{}}
        \toprule
        Method
        & $f^\star$
        & $p$
        & $c$
        & $L$
        & $U_0$
        & Stepsize schedule
        & Reference \\
        \midrule
        \DCSG{}
        & --
        & $\checkmark$
        & $\checkmark$
        & --
        & $\checkmark$
        & Adaptive
        & This work \\

        Polyak
        & $\checkmark$
        & --
        & --
        & --
        & --
        & Adaptive
        & \cite{Polyak1969} \\

        RSG
        & --
        & $\checkmark$
        & $\checkmark$
        & $\checkmark$
        & $\checkmark$
        & Stagewise
        & \cite{YangLin2018} \\

        DS-SG
        & --
        & $\checkmark$
        & $\checkmark$
        & $\checkmark$
        & $\checkmark$
        & Stagewise
        & \cite{JohnstoneMoulin2020} \\

        DecayGD
        & --
        & $\checkmark$
        & $\checkmark$
        & $\checkmark$
        & $\checkmark$
        & Predetermined
        & \cite{Wang2025} \\
        \bottomrule
    \end{tabular}
\end{table}

\subsection{Power objectives with exact growth constants} \label{subsection:power-objective}
We first consider a family of power objectives with an explicit minimizer
and exact growth constants: 
\begin{equation} \label{eq:experiment-1}
    \min_{x\in C} f_p(x)
    := \frac{1}{p}\|A(x-x^\star)\|_1^p,
    \qquad
    C:=\{x\in\mathbb{R}^d:\|x\|_2\leq 1\},    
\end{equation}
where $d=50$ and $p\in\{1,2,3\}$. We construct $A=Q_1\operatorname{diag}(\sigma)Q_2^\top$, where $Q_1$ and $Q_2$
are random orthogonal matrices and the singular values are logarithmically spaced between $1/2$ and $1$. Thus, $\|A\|_2=1$ and $\operatorname{cond}_2(A)=2$. We choose $x^\star$ such that $\|x^\star\|_2=1/2$. By construction, $A$ is nonsingular, and thus the solution set is $X^\star=\{x^\star\}$ and the optimal value is $f^\star=0$.

The objective in \cref{eq:experiment-1} satisfies the H\"older growth condition \cref{eq:holderian-error-bound} with 
$
    c=\frac{1}{p}    \left(\max_{1\leq j\leq d}\|A^{-1}e_j\|_2\right)^{-p},
$ 
where $e_j$ denotes the $j$th standard basis vector. We use the subgradient
$$
    g(x)=\|A(x-x^\star)\|_1^{p-1}
    A^\top\operatorname{sign}\bigl(A(x-x^\star)\bigr),
$$ 
with $\operatorname{sign}(0)=0$. This oracle satisfies the shrinking-subgradient condition \cref{eq:shrinking-selected-subgradients} with $s=p-1$. A uniform subgradient bound on $C$ is $L=2^{p-1}d^{p/2}$.
All methods are initialized at $x_0=0$, and we use the initial squared-distance bound $U_0=1$. This bound is conservative, since $e_0=\|x_0-x^\star\|_2^2=1/4$. We run each method for $N=4000$ subgradient evaluations. 

\begin{figure}[t]
\setlength{\abovecaptionskip}{2pt}
  \centering
  \includegraphics[width = \textwidth]{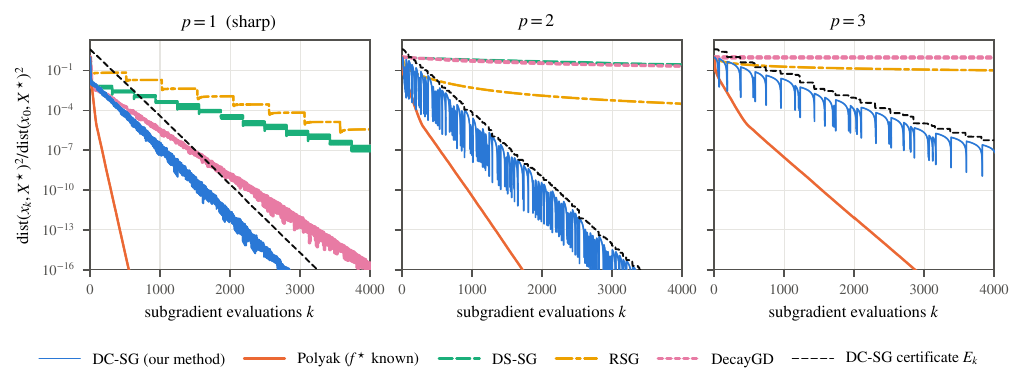}
  \caption{Problem \cref{eq:experiment-1} with $p\in\{1,2,3\}$: Normalized squared distance $e_k/e_0$ versus the number of subgradient evaluations for \DCSG{} (\Cref{alg:dc-sg}), Polyak method, DS-SG, RSG, and DecayGD, where $e_k=\operatorname{dist}(x_k,X^\star)^2$. The black dashed curve shows the normalized \DCSG{} certificate $E_k/e_0$.}
  \label{fig:exp1_power}
\end{figure}

The numerical results are shown in \cref{fig:exp1_power}, which are consistent with the theoretical guarantees. In all three cases, the certificate satisfies $e_k\leq E_k$ throughout the iterations of \DCSG{} and decreases geometrically. For $p=1$, this behavior agrees with \Cref{thm:geometric-certificate-decay-p-equals-1}; for $p=2,3$, it agrees with~\Cref{cor:shrinking-subgradient-rates} under the
shrinking condition $s=p-1$. As shown in \cref{fig:exp1_power}, although the squared-distance error of \DCSG{} may oscillate, its certificate gives a monotonically decreasing upper bound, as guaranteed by \cref{prop:certificate-monotonicity}. 
In this experiment, the Polyak method using with the exact optimal value $f^\star=0$ converges fastest among all methods. Our \DCSG{} also exhibits linear convergence for all three values of $p$ and maintains a valid distance certificate without using $f^\star$. In addition, \DCSG{}
converges much faster than DecayGD, DS-SG, and RSG. In particular, the methods using predetermined schedules (DecayGD, DS-SG, and RSG) make substantially less progress for $p=2,3$. This comparison clearly shows the benefit of using the observed subgradient norm in the stepsize selection.

\subsection{Regularized support vector machines}
We next consider an $\ell_2$-regularized support vector machine with hinge loss:
\begin{equation}\label{eq:experiment-2}
    \min_{w\in C} f(w)
    := \frac{1}{m}\sum_{i=1}^{m}\max\{0,1-y_i a_i^\top w\}
       +\frac{\lambda}{2}\|w\|_2^2,
    \qquad
    C:=\{w\in\mathbb{R}^{d}:\|w\|_2\leq R\},
\end{equation}
where $\lambda=10^{-2}$ and $R=\sqrt{2/\lambda}$. We use the Wisconsin Diagnostic Breast Cancer
dataset~\cite{breastCancerWisconsin}, which contains $m=569$ samples with $30$ features. We standardize each feature to have zero mean and unit variance, and append a constant coordinate to include an intercept.

\begin{figure}[t]
\setlength{\abovecaptionskip}{0pt}
  \centering
  \includegraphics[width = 0.78\textwidth]{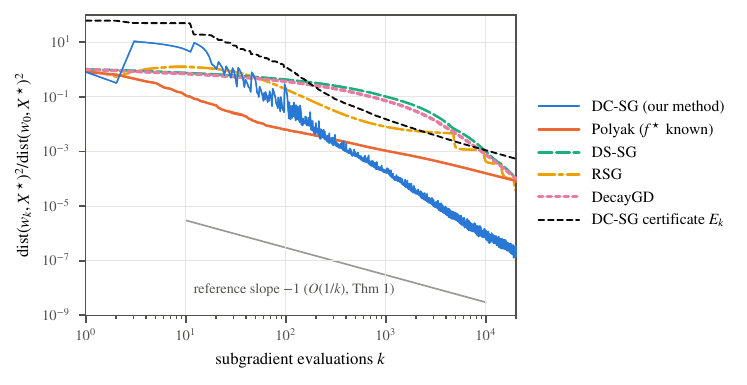}
  \caption{Problem \cref{eq:experiment-2} with  $\lambda=10^{-2}$: normalized squared distance $e_k/e_0$
    versus the number of subgradient evaluations, where  $e_k=\|w_k-w^\star\|_2^2$. We compare \DCSG{}, the Polyak method, DS-SG, RSG, and DecayGD. The black dashed curve shows the normalized \DCSG{} certificate $E_k/e_0$.}
  \label{fig:exp2_svm}
  \vspace{-1mm}
\end{figure}

The objective in \cref{eq:experiment-2} is $\lambda$-strongly convex and has a unique minimizer $w^\star$.
Moreover, $\lambda\|w^\star\|_2^2/2\leq f(w^\star)\leq f(0)=1$, so $w^\star\in C$. Consequently, the H\"older growth condition \cref{eq:holderian-error-bound} holds with~$p=2$ and $c=\lambda/2=5\times10^{-3}$. We use the subgradient
$
    g(w)=-\frac{1}{m}\sum_{i:y_i a_i^\top w<1}y_i a_i
          +\lambda w. 
$ The selected subgradients need not vanish as $w$ approaches
$w^\star$. Applying \Cref{thm:polynomial-certificate-decay-lipschitz} with $p=2$ yields $E_k=O(1/k)$ for \DCSG.  We compute a reference
solution using an interior-point solver, Clarabel \cite{goulart2026clarabel}.
The reference solution is used to evaluate the errors, and its objective value is also used in the Polyak method. All methods are initialized at $w_0=0$, and we use the initial
squared-distance bound $U_0=R^2=200$. We run each method for $N=20{,}000$ iterations.

The numerical results are shown in \cref{fig:exp2_svm}. The distance certificate from \DCSG{} exhibits an approximately $\mathcal{O}(1/k)$ rate over the later iterations. This behavior is consistent with \Cref{thm:polynomial-certificate-decay-lipschitz} for $p=2$. The measured squared-distance error decreases faster than the certificate, indicating that the certificate provides a conservative estimate of the actual error in this example. After an initial transient, \DCSG{} achieves substantially higher distance accuracy than the comparison methods. At $N=20{,}000$, its normalized squared distance is $1.5\times10^{-7}$, while those of the other methods range from $3.0\times10^{-5}$ to $9.1\times10^{-5}$.
This illustrates the effectiveness of \DCSG{} on a nonsmooth problem with an explicit growth constant and uniformly bounded subgradients. 

\section{Conclusion} \label{section:conclusion}
\vspace{-1mm}
In this paper, we have proposed \DCSG{}, a new distance-certified subgradient method
for convex optimization under H\"older growth. Our \DCSG{} method maintains a computable upper bound on the squared distance to the solution set and uses this certificate to select Polyak-type stepsizes, without requiring the optimal objective value. We have established sublinear convergence for $p>1$ and linear
convergence for $p=1$, together with sharper rates under shrinking subgradients, including linear convergence when $s=p-1$. The numerical experiments support these theoretical findings. For future work, an interesting direction is to develop adaptive estimates of the growth parameters while preserving the validity of the distance certificate. 
\vspace{5mm}

\noindent \textbf{Statements on AI use.}
This paper started with a careful examination of the faster subgradient methods in \cite{JohnstoneMoulin2020} in March 2026 as an SRIP topic for the first author. Motivated by the Polyak stepsize, the authors had the idea of the one-step minimax principle for choosing the stepsize in subgradient methods under H\"older growth. We used ChatGPT 5.6 and 6 to assist with some~derivations, literature searches, proof checking, and language editing. We also used Claude Opus 5 for coding and numerical experiments. The authors carefully reviewed all AI-assisted material and independently verified the mathematical arguments. The authors assume responsibility for all content.

\bibliographystyle{unsrtnat}
\bibliography{reference}

\end{document}